\documentclass[oneside,11pt,a4paper]{article}

\usepackage[round,sort&compress]{natbib}

\usepackage[margin=1in]{geometry}

\usepackage{url}
\usepackage[colorlinks,bookmarksopen,bookmarksnumbered,citecolor=blue,linkcolor=blue,urlcolor=blue]{hyperref}

\usepackage{multicol}
\usepackage{graphicx}
\usepackage{subcaption}
\usepackage{amssymb,amsfonts}
\usepackage{amsmath}
\usepackage{mathrsfs}

\usepackage{yfonts}
\usepackage{changepage}

\usepackage{authblk}

\title{Bayesian statistical learning using density operators}

\author{Yann Berquin\thanks{\href{mailto:yann.berquin@ccnu.edu.cn}{yann.berquin@ccnu.edu.cn}}}
\affil{CCNU Wollongong Joint Institute\\Central China Normal University\\Wuhan, China}
       
\date{}

\begin{document}

\maketitle

\begin{abstract}
This short study reformulates the statistical Bayesian learning problem using a quantum mechanics framework. Density operators representing ensembles of pure states of sample wave functions are used in place probability densities. We show that such representation allows to formulate the statistical Bayesian learning problem in different coordinate systems on the sample space. We further show that such representation allows to learn projections of density operators using a kernel trick. In particular, the study highlights that decomposing wave functions rather than probability densities, as it is done in kernel embedding, allows to preserve the nature of probability operators. Results are illustrated with a simple example using discrete orthogonal wavelet transform of density operators.
\end{abstract}

	\section{Introduction}\label{sec:I1}
	
 Consider an unknown probability distribution on a sample space $\mathscr{S}$ and a set of $N$ independent noisy singletons $\{S_i\}_{i \in I}$ generated by this unknown probability distribution. We interest ourselves in the problem of recovering the unknown probability distribution from the space $\mathscr{P}$ of all absolutely continuous probability distributions given this set of independent noisy singletons. $\mathscr{P}$ can be identified without ambiguity to the set of all probability densities on $\mathscr{S}$ as densities are almost unique for a given probability distribution. Singleton noises are assumed to be independent from probability distributions such that $P( S_i | Z, s)=P( S_i | s)$ for all $Z \in \mathscr{P}$. The resulting probability distribution $P( S_i | s)$ represents the probability of drawing the singleton (or sample) $S_i$ with an unknown true value $s$. If $\mathscr{P}$ is measurable, it is then possible to define probability distributions on $\mathscr{P}$ and apply Bayes' rule. Given a prior probability distribution $P(Z)$ on $\mathscr{P}$, the corresponding Bayesian statistical learning problem reduces to finding the posterior probability distribution $P \big( Z| \{S_i\}_{i \in I} \big)$ on $\mathscr{P}$ where
\begin{equation}
\label{eq:learning0}
\begin{aligned}
&P \big( Z| \{S_i\}_{i \in I} \big) = \frac{P( Z )}{P \big(\{S_i\}_{i \in I} \big)} \, \prod_{i \in I} \int_{\mathscr{S}} P \big( S_i | s \big) \, Z(s)
\end{aligned}
\end{equation}
Equation \ref{eq:learning0} provides the general solution to the probability distribution learning problem with independent samples. In essence, it computes the joint likelihood to generate independent samples $\{S_i\}_{i \in I}$ from the probability distribution $Z$ in conjunction with a prior expectation over $\mathscr{P}$. In this study, we propose to replace probability distributions with density operators which implies that singletons can be described by wave functions. The idea is far from being novels and numerous studies have already tried to use quantum mechanics framework to describe classical statistical problems. For instance in the work by \citet{wolf2006learning}, density operators are used to perform spectral clustering. Quantum mechanics concepts have also been applied to specific machine learning algorithms \citep[e.g.][]{gonzalez2021classification}. In particular, quantum mechanics concepts have been used for nearest mean classifiers \citep[e.g.][]{sergioli2018quantum}, binary classifiers \citep[e.g.][]{tiwari2019binary} and sentiment learnings \citep[e.g.][]{zhang2018quantum,li2021quantum}. The connection between quantum mechanics and classical probabilities has also been largely investigated \citep[e.g.][]{malley1993quantum,barndorff2003quantum,jarzyna2020geometric} but is beyond the scope of our work. In this study, quantum mechanics is merely used to decompose samples on a complete orthonormal basis of the sample space. The novelty of this preliminary study lies in the reformulation of the Bayesian statistical learning using density operators in different basis (i.e. coordinates) and its application to learning projection of density operators. In particular, the study highlights that decomposing wave functions rather than probability densities, as is done in kernel embedding, allows to preserve the nature of density operators. Section \ref{sec:2} discusses mostly the mathematical background and recalls the link between probability distributions and density operators. Section \ref{sec:3} introduces the general Bayesian statistical learning of density operators and details how it can be used to learn projection of density operators. Results are illustrated with discrete wavelet transforms of density operators.

	\section{Mathematical background}\label{sec:2}
	
	 This section recalls fundamentals of quantum mechanics and shows how density operators can be used in place of probability distributions and measurement operators can be used to represent non-noisy singletons. For the sake of clarity, the sample space $\mathscr{S}$ is first assumed to be a discrete collection of points such that $\mathscr{S} = \{s_j\}_{j \in K}$ where the size of $K$ is a set of integers which can be infinitely large. The space of all square integrable complex functions on the sample space $\mathscr{S}$ forms a closed linear subspace of the sequence space $\ell^2$. It is thus a Hilbert space which shall be denoted $\mathscr{H}$. The generalization to continuous sample spaces is discussed in section \ref{sec:2-2}
	 
	 \subsection{Probability distributions as density operators}
	 
In quantum mechanics, wave functions are probability amplitudes which are defined as complex functions on $\mathscr{S}$ which modulus square are probability densities. Wave functions can be identified to unit vectors of the Hilbert space $\mathscr{H}$. Using the set of samples $\{s_j\}_{j \in J}$ it is possible to define the family of wave functions $\{|s_j\rangle\}_{j \in J}$ such that
\begin{equation}
\begin{aligned}
\label{eq:def_01}
&\langle s_j | s_k \rangle = \delta_{jk}
\end{aligned}
\end{equation}
where $\delta_{jk}$ is the usual Kronecker delta function. By construction, wave functions $\{|s_j\rangle\}_{j \in J}$ are mutually orthonormal and form a complete basis of $\mathscr{H}$. Using this family of wave functions, any wave functions $|z\rangle$ associated to a given probability distribution $Z$ can be written as follows
\begin{equation}
\begin{aligned}
\label{eq:def_02}
&|z\rangle = \sum_{j \in J}\langle s|z \rangle \, |s\rangle
\end{aligned}
\end{equation}
where $\{\langle s_j|z \rangle\}_{j \in J}$ is a set of complex coefficients such that $|\langle s_j|z \rangle|^2=Z(s_j)$ for all $j \in J$. There exists an infinite number of wave functions associated to a given probability density as each complex coefficient $\langle s_j|z \rangle$ is only constrained by its squared modulus. Wave functions are essential bricks of quantum mechanics along with density operators and positive-operator-valued measures (POVM). Density operators are positive semi-definite, Hermitian operators of trace one which provide a description of quantum systems. Given the family of mutually orthonormal wave functions $\{|s_j\rangle\}_{j \in J}$, any density operator $\rho$ acting on $\mathscr{H}$ can be expressed as
\begin{equation}
\begin{aligned}
\label{eq:def_03}
&\rho =  \sum_{l \in J} \sum_{j \in J}  w(s_j,s_l) \,  | s_j\rangle \langle s_l| 
\end{aligned}
\end{equation}
We shall denote $\mathscr{G}$ the set of all density operators acting on $\mathscr{H}$. Since $\rho$ is positive semi-definite, it is easy to verify that $w(s_j,s_l) = \overline{w(s_j,s_l)}$ for any $j \in J$ and $l \in J$. In addition, since $\mathrm{tr}(\rho)=1$, it is also easy to verify that $\sum_{j \in J} w(s_j,s_j) =1$. In quantum mechanics, measurements on quantum systems (i.e. density operators) are represented by POVM \citep[e.g.][]{paris2012modern}. Loosely speaking, a POVM is a set of positive semi-definite operators acting on wave functions that sum to the identity matrix. In particular, the set of operator $\{|s_j \rangle \langle s_j |\}_{j \in J}$ associated to each sample on the sample space corresponds to such a measurement operator. Since the sample space is discrete, $\{|s_j \rangle \langle s_j |\}_{j \in J}$ corresponds to a projection-valued measure and is a projection operator. Interestingly, the projection operator associated to the measurement of a sample in state $|s_j\rangle$ is also a density operator. It can be interpreted as the density operator of a single particle (i.e. non-noisy singleton) in state $|s_j \rangle$. The probability of measuring a sample in state $|s_j\rangle$ given density operator $\rho$ can be obtained using Born's rule
\begin{equation}
\begin{aligned}
\label{eq:def_06}
&P (s_j |\rho ) = \mathrm{tr} \Big( \rho \, |s_j \rangle \langle s_j | \Big) =  w(s_j,s_j) 
\end{aligned}
\end{equation}
Note that we (abusively) use the notation $P (s_j |\rho )$ instead of $P \Big(|s_j \rangle|\rho \Big)$ to make it easier to read. Consider now a family of mutually orthonormal wave functions $\{|\psi_j\rangle\}_{j \in J}$ which also forms a complete orthonormal basis of $\mathscr{H}$. Using this family of wave functions, any density operator $\rho$ can be expressed as
\begin{equation}
\begin{aligned}
\label{eq:def_07}
&\rho =  \sum_{j \in J} \sum_{l \in J} w(\psi_j ,\psi_l)  \, |\psi_j \rangle \langle \psi_l |
\end{aligned}
\end{equation}
where
\begin{equation*}
\begin{aligned}
w(\psi_j ,\psi_l) = \sum_{k \in J} \sum_{n \in J} w(s_k,s_n)  \, \langle \psi_j|s_k\rangle \, \langle s_n | \psi_l \rangle
\end{aligned}
\end{equation*}
Since $\{|\psi_j\rangle\}_{j \in J}$ are mutually orthonormal wave functions, coefficients $\{w(\psi_j ,\psi_l)\}_{j,l}$ uniquely determine density operator $\rho$. It is also easy to verify that $w(\psi_j ,\psi_l)= \overline{w(\psi_j ,\psi_l)}$ for any $j \in J$ and $l \in J$ and $\sum_{j \in J} w(\psi_j ,\psi_j) =1$. Using the decomposition of $\rho$ on $\{| \psi_j \rangle\}_{j \in J}$, probability $P(s_j |\rho)$ can be expressed as a function of coefficients $\{w(\psi_j ,\psi_l)\}_{j,l}$
\begin{equation}
\begin{aligned}
\label{eq:def_08}
P(s_j |\rho) &= p \Big( s_j | \{w(\psi_j ,\psi_l)\}_{j,l} \Big) \\
&=\sum_{k \in J} \sum_{l \in J} w(\psi_k ,\psi_l) \, \langle s_j |\psi_k \rangle \, \langle \psi_l | s_j \rangle
\end{aligned}
\end{equation}
For a given probability distribution $Z$, any density operator $\rho$ which satisfies $P (s_j|\rho )=Z(s_j)$ for all $j \in J$ can be associated to the probability distribution $Z$. There exists an infinite number of density operators associated to a given probability density. This can be easily verified since for any wave function $|z\rangle$ associated to a given probability distribution $Z$, $|z\rangle \langle z|$ is a density operator associated to $Z$. For a given probability distribution, of notable interest is the density operator such that $w(s_j,s_l) = 0$ if $j\neq l$ for all $j$ and $l$ in $J$. Such density operator corresponds to a quantum system composed of an ensemble of samples $\{|s_j\rangle\}_{j \in J}$ occurring with probability $Z(s_j)$ and can be written as
\begin{equation}
\begin{aligned}
\label{eq:def_09}
\rho &= \sum_{j \in J} Z(s_j) \, | s_j \rangle \langle s_j |
\end{aligned}
\end{equation}
The set of all such density operators will be denoted by $\mathscr{E}$ such that $\mathscr{E} \in \mathscr{G}$ where $\mathscr{G}$ denotes the set of all density operators acting on $\mathscr{H}$.
\begin{figure}
    \centering
        \includegraphics[width=0.6\textwidth]{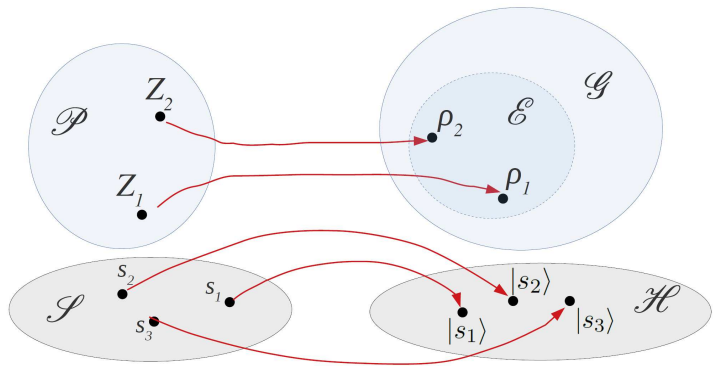}
    \caption{Correspondence between classic probability elements and their quantum analogs. In this study, we constrain density operators to the subset $\mathscr{E}$ of all density operators $\mathscr{G}$ acting on $\mathscr{H}$.}\label{fig:map}
\end{figure}
Representing probability distributions as quantum systems composed of an ensemble of pure states $\{|s_j\rangle\}_{j \in J}$ ensures a one to one correspondence between density operators and probability distributions (see Figure \ref{fig:map}). In the rest of this study, we shall limit ourselves to density operators corresponding to quantum systems composed of an ensemble of pure states $\{| s_j \rangle\}_{j \in J}$.

\subsection{Generalizing to continuous sample spaces}\label{sec:2-2}

The sample space $\mathscr{S} = \{s_j\}_{j \in J}$ has been assumed discrete so far which is very restrictive in practice. It is possible to extend the validity of the study to continuous sample spaces for which the space of all square integrable complex functions on the sample space $\mathscr{S}$ forms the Hilbert space $L^2(\mathscr{S},\mu)$ where $\mu$ is a strictly positive finite Borel measure. Using the set of samples on $\mathscr{S}$ it is possible to define the family of wave functions $\{|s\rangle\}_{s \in \mathscr{S}}$ such that
\begin{equation}
\begin{aligned}
\label{eq:def_11}
&\langle s | s' \rangle = \delta(s-s') \quad \mathrm{and} \quad \langle s | \phi \rangle = \phi(s)
\end{aligned}
\end{equation}
where $\delta(s)$ is the Dirac delta function and $| \phi \rangle$ a wave function associated to any function $\phi \in L^2(\mathscr{S},\mu)$, i.e. $\int_\mathscr{S} |\phi(s)|^2 \, d\mu(s)=1$. Note that the inner product now corresponds to the usual inner product on $L^2(\mathscr{S},\mu)$ such that for any wave functions $|\phi_1\rangle$ and $|\phi_2\rangle$ the inner product can be explicitly written as $\langle \phi_1|\phi_2\rangle = \int_\mathscr{S} \overline{\phi_1(s)} \, \phi_2(s) \, d\mu(s)$. By construction, wave functions $\{|s\rangle\}_{s \in \mathscr{S}}$ are mutually orthogonal and correspond to the well-known eigenvector of the position operator \citep[e.g.][]{phillips2013introduction}. According to equation \ref{eq:def_11} wave functions $\{|s\rangle\}_\mathscr{S}$ are such that $|\langle s|s \rangle| = \delta(0)$ and hence do not seem to correspond to unit vectors of $L^2(\mathscr{S},\mu)$. This issue is well-known and is related to the definition of the Dirac delta function in \ref{eq:def_11}. We propose in this study to use the hyperreal delta function introduced in \citet{cabbolet2021hyperreal} such that $\delta(0)= \omega$ where $\omega$ is the usual positive infinitely large hyperreal number with $|\omega| = \infty$. This definition allows to introduce the following hyperreal density operator $\rho$ which generalizes equation \ref{eq:def_09} to continuous sample space $\mathscr{S}$
\begin{equation}
\begin{aligned}
\label{eq:def_12}
\rho = \frac{1}{\omega} \, \int_{\mathscr{S}} Z \left(s\right) \, | s\rangle \langle s |
\end{aligned}
\end{equation}
It is straightforward to verify that this definition yields $\mathrm{tr}  (\rho)=1$. In addition, using measurement operator $|s \rangle \langle s |$, it can be found that the probability density of measuring a sample in state $|s\rangle$ given density operator $\rho$ in equation \ref{eq:def_12} is
\begin{equation}
\begin{aligned}
\label{eq:def_13}
p(s |\rho) = \mathrm{tr}  \Big( |s \rangle \langle s| \, \rho \Big) = \zeta \left( s\right)
\end{aligned}
\end{equation}
where $\zeta$ is the usual probability density (i.e. Radon–Nikodym derivative) associated to the probability distribution $Z$ and measure $\mu$. Note that while $|s \rangle \langle s|$ corresponds to a projection-valued measure when the sample space is discrete, in the case of a continuous sample space $|s \rangle \langle s|$ corresponds to a positive operator-valued measure. As a result, when the sample space is continuous, $|s \rangle \langle s|$ is not a projection operator anymore. Table \ref{table:definition} provides the correspondence between classic probability elements and associated quantum analogs according to the conventions chosen. There is a strict one-to-one correspondence between probability densities and density operators describing quantum systems composed of ensemble of pure states $\{|s\rangle\}_{s \in \mathscr{S}}$. There is also a strict one-to-one correspondence between singletons $\{S_i\}_{i \in I}$ (i.e. points on the sample space) and wave functions $\{|S_i\rangle\}_{i \in I}$ (see Figure \ref{fig:map}).
\begin{table}[ht]
\centering
\caption{Classic formalism vs quantum analog.}
\renewcommand{\arraystretch}{1.2}
\begin{tabular}{ l l l }
  \hline
   & Classic  & Quantum  \\
   \hline
  Representing probability & $Z(s)$ & $\frac{1}{\omega} \int_\mathscr{S} Z(s) |s\rangle \langle s|$ \\
  Representing samples & $\{S_i\}_{i \in I}$ & $\{|S_i\rangle\}_{i \in I}$ \\
  \hline
\end{tabular}
\label{table:definition}
\end{table}
Consider now a family of mutually orthonormal wave functions $\{|\psi_j\rangle\}_{j \in J}$ which forms a complete orthonormal basis of $L^2(\mathscr{S},\mu)$. Using this family of wave functions, the density operator $\rho$ can be expressed as
\begin{equation}
\begin{aligned}
\label{eq:def_14}
&\rho = \frac{1}{\omega} \, \sum_{j \in J} \sum_{l \in J} w(\psi_j,\psi_l) \, |\psi_j \rangle \langle \psi_l |
\end{aligned}
\end{equation}
where
\begin{equation*}
\begin{aligned}
\label{eq:def_14_1}
w(\psi_j,\psi_l) &= \int_{\mathscr{S}} Z(s) \, \langle s |\psi_l \rangle \, \langle \psi_j |s \rangle
\end{aligned}
\end{equation*}
Density operator $\rho$ as defined in equation \ref{eq:def_12} allows to represent any probability distribution $Z$ as a quantum system composed of an ensemble of pure states $\{|s\rangle\}_{s \in \mathscr{S}}$. For a given complete orthonormal basis $\{|\psi_j\rangle\}_{j \in J}$ of $L^2(\mathscr{S},\mu)$ and $\rho \in \mathscr{E}$, the set of complex coefficients $\{w(\psi_j,\psi_l)\}_{j,l}$ as defined in equation \ref{eq:def_14_1} can be understood as coordinates of the probability distribution $Z$.\newline
\\
We have detailed how probability distributions can be represented using density operators as quantum systems composed of an ensemble of pure states $\{|s\rangle\}_{s \in \mathscr{S}}$. Using this formalism, the main difference with the classical approach is the nature of the sample space where samples are treated as wave functions. Describing samples as wave functions allows to naturally represent density operators in any orthonormal basis $\{|\psi_j\rangle\}_{j \in J}$. In the following section we detail how such representation can be used to formulate the associated Bayesian statistical learning problem.

	\section{Statistical learning of density operators}\label{sec:3}

We interest ourselves in finding the density operator in $\mathscr{E}$ associated to an unknown probability density, in any given basis, given a set of independent noisy singletons $\{S_i\}_{i \in I}$ generated by the unknown probability density. This merely corresponds to the usual statistical learning problem where density operators in $\mathscr{E}$ are used in place of probability distributions. Bayes' rule allows to us to write
\begin{equation}
\begin{aligned}
\label{eq:res_0}
&P \left(\rho \big| \{S_i\}_{i \in I}  \right) \propto P \left( \rho\right) \, \prod_{i \in I} \mathrm{tr} \Big( M(S_i) \, \rho \Big)
\end{aligned}
\end{equation}
where $M(S_i)$ can be understood as the measurement operator associated to each singleton $S_i$ endowed with its associated measurement uncertainties $P(S_i|s)$ such that
\begin{equation}
\begin{aligned}
M(S_i) &= \int_{\mathscr{S}} P(S_i|s) \, |s\rangle \langle s|\\
\mathrm{tr} \Big( M(S_i) \, \rho \Big) &= \int_{\mathscr{S}} P\left(S_i|s\right) \,  p \left(s \big| \rho \right) \, d\mu(s)
\end{aligned}
\end{equation}
Since there is a strict one-to-one between probability densities and density operators in $\mathscr{E}$, equation \ref{eq:res_0} is fully equivalent to equation \ref{eq:learning0} where probability density $\zeta$ has been replaced with its associated density operator $\rho$. $P(S_i|s_i)$ is the usual measurement uncertainty associated to $S_i$, $P( \rho )=P( \zeta)$ is the usual prior probability distribution and $p \left(s \big| \rho \right)=\zeta(s)$. The interesting property arises from the possibility to naturally express density operators in any basis of $L^2(\mathscr{S},\mu)$. In particular, for any given complete orthonormal basis $\{|\psi_j\rangle\}_{j \in J}$ of $L^2(\mathscr{S},\mu)$ with the set of complex coefficients $\{w(\psi_j,\psi_l)\}_{j,l}$ associated to $\rho \in \mathscr{E}$, probabilities in equation \ref{eq:res_0} reduce to functions of the complex coefficients such that $P (\{w(\psi_j,\psi_l)\}_{j,l} | \{S_i\}_{i \in I}  )=P( \rho | \{S_i\}_{i \in I}  )$, $P( \{w(\psi_j,\psi_l)\}_{j,l} )=P \left( \rho\right)$ and $p \left(s \big| \{w(\psi_j,\psi_l)\}_{j,l} \right)=p \left(s \big| \rho \right)$. As a result, Bayesian statistical learning of density operators introduced in equation \ref{eq:res_0} can be also written as
\begin{equation}
\label{eq:res_1}
\boxed{\boxed{P \left(\{w(\psi_j,\psi_l)\}_{j,l} \big| \{S_i\}_{i \in I}  \right) \propto P \left( \{w(\psi_j,\psi_l)\}_{j,l} \right) \, \prod_{i \in I} \int_{\mathscr{S}} P\left(S_i|s\right) \, p \left(s \big| \{w(\psi_j,\psi_l)\}_{j,l} \right) \, d\mu(s)}}
\end{equation}
The usual statistical learning problem as described in equation \ref{eq:learning0} thus appears as nothing more than a special case of the statistical problem described in equation \ref{eq:res_1} with the complete orthonormal basis chosen as the ensemble of pure states $\{|s\rangle\}_{s \in \mathscr{S}}$.\newline
\\
\textbf{Homogeneous prior with non-noisy singleton.} \textit{We provide an expression for the mode of the posterior probability distribution $P (\{w(\psi_j,\psi_l)\}_{j,l} | \{S_i\}_{i \in I} )$ when the prior probability distribution $P ( \{w(\psi_j,\psi_l)\}_{j,l} )$ is homogeneous (i.e. constant) and when $p(S_i|s) = \delta(S_i-s)$ for all singletons. The mode $\tilde{\rho}$ of the posterior probability distribution $P (\{w(\psi_j,\psi_l)\}_{j,l} | \{S_i\}_{i \in I} )$ can then be expressed using equation \ref{eq:def_14} as
\begin{equation}
\begin{aligned}
\label{eq:res_2}
\tilde{\rho} = \frac{1}{\omega} \, \sum_{j \in J} \sum_{l \in J} \tilde{w}(\psi_j,\psi_l) \, |\psi_j \rangle \langle \psi_l |
\end{aligned}
\end{equation}
where
\begin{equation*}
\begin{aligned}
w(\psi_j,\psi_l) &= \int_{\mathscr{S}} \tilde{\zeta}(s) \, \langle s |\psi_l \rangle \, \langle \psi_j |s \rangle \, d\mu(s)\\
&= \frac{1}{N} \, \sum_{i \in I} \langle S_i |\psi_l \rangle \, \langle \psi_j |S_i \rangle\\
&= \langle \psi_j | \left( \frac{1}{N} \, \sum_{i \in I} |S_i \rangle \langle S_i| \right) |\psi_l \rangle 
\end{aligned}
\end{equation*}
$\tilde{\zeta}$ denotes the mode of $P(\zeta|\{S_i\}_{i \in I})$ (see equation \ref{eq:learning0}) such that $\tilde{\zeta}(s)= 1/N \, \sum_{i \in I} \delta(S_i-s)$. We recall that $N$ corresponds to the number of singletons in $\{S_i\}_{i \in I}$ (i.e. $N=|I|$).}\newline
\\
Using wave functions to represent samples provides a natural way to decompose samples on any orthonormal basis. Furthermore, since density operators are described as ensembles of wave functions, density operators can also be decomposed on any complete wave functions basis.

	\subsection{Learning with embedded samples}\label{sec:mapping}

Consider a family $\{|\psi_j \rangle \}_{j \in J}$ of mutually orthonormal wave functions which forms a complete orthonormal basis of $L^2(\mathscr{S},\mu)$ and a set of complex coefficients $\{\alpha_j\}_{j \in J}$. Using the subset of wave functions $\{|\psi_j \rangle \}_{j \in L}$ allows to introduce the following linear operator $A$
\begin{equation}
\begin{aligned}
\label{eq:ker_0}
A = \sum_{j \in L} \alpha_j \, |\psi_j \rangle \langle \psi_j |
\end{aligned}
\end{equation}
Such linear operator can be loosely understood as a quantum operator which maps (up to a normalization coefficient) any wave function $\sum_{j \in J} c_j \, |\psi_j\rangle$ to a new wave function $\sum_{j \in L} \alpha_j \, c_j \, |\psi_j\rangle$. In particular, to each wave function $|s\rangle$, it associates a new wave function $A | s \rangle$. Note that the set of wave functions $\{A| s\rangle\}_{s \in \mathscr{S}}$ are not necessarily mutually orthogonal. For each density operator $\rho \in \mathscr{E}$ which are not in the kernel of $A$, we let $\rho_A$ be such that,
\begin{equation}
\begin{aligned}
\label{eq:ker_1}
\rho_A = \frac{A \, \rho \, A^*}{\mathrm{tr}\left( A \, \rho \, A^* \right)}
\end{aligned}
\end{equation}
where $A^*$ denotes the conjugate transpose of $A$. The set of all density operators $\rho_A$ forms a subspace of all density operators operating on $L^2(\mathscr{S},\mu)$ which is not necessarily in $\mathscr{E}$. Interestingly, the operator $K=A^*A=A A^*$  is a symmetric, positive definite operator with eigenfunctions $\{|\psi_j \rangle \}_{j \in J}$ and eigenvalues $\{|\alpha_j|^2 \rangle \}_{j \in J}$. The operator $K$ can thus be associated to a reproducing kernel Hilbert space (RKHS) $H$ such that for any two wave functions $\phi =\sum_{j \in L}  c_{j} \, |\psi_j\rangle$ and $\phi' =\sum_{j \in L}  c_{j}' \, |\psi_j\rangle$ it is possible to write without ambiguity $\langle \phi | \phi' \rangle_H = \langle \phi | K | \phi' \rangle$. As a result, for the density operator associated to the probability distribution $Z$ which corresponds to a quantum system composed of an ensemble of pure states $\{|s\rangle\}_{s \in \mathscr{S}}$ (see equation \ref{eq:def_12}), the density operator $\rho_A$ is such that
\begin{equation}
\begin{aligned}
\label{eq:ker_2}
\rho_A &= \frac{1}{\mathrm{tr}\left( A \, \rho \, A^* \right)} \,  \int_{\mathscr{S}}  Z \left(s\right) \,A| s\rangle \langle s | A^*\\
&= \frac{1}{\mathrm{tr}\left( A \, \rho \, A^* \right)} \, \int_{\mathscr{S}}  \Big(  Z \left(s\right) \, \langle s |K| s\rangle \Big) \,\widehat{A| s\rangle} \widehat{\langle s | A^*}
\end{aligned}
\end{equation}
where
\begin{equation*}
\begin{aligned}
\widehat{A| s\rangle} = \frac{A| s\rangle}{\sqrt{\langle s | K| s\rangle}}
\end{aligned}
\end{equation*}
$\widehat{A| s\rangle}$ denotes the normalized wave function associated to $A|s\rangle$ which can now be interpreted as the embedding of $|s\rangle$ on $H$. According to equation \ref{eq:ker_2}, density operator $\rho_{A}$ can be understood as an ensemble of embedded samples $\widehat{A| s\rangle}$ with probability $Z \left(s\right) \, \langle s |K| s\rangle$. It is worth noting that when the kernel is such that $\langle s |K| s\rangle$ is constant for all $s \in \mathscr{S}$, then the probability of embedded samples $\widehat{A| s\rangle}$ in $\rho_A$ is the same as the probability of samples $|s\rangle$ in $\rho$.\newline
\\
We now wish to find the probability density $p( A| s\rangle  |\rho_{A})$ that the system described by $\rho_A$ yields measurement in state $A| s\rangle$. The set of bounded non-negative self-adjoint operators $\{A |s\rangle  \langle s| A^*\}_{s \in \mathscr{S}}$ is such that $\int_\mathscr{S} A |s\rangle  \langle s| A^* =K$. Since $K$ is the reproducing kernel of $H$, $\{A |s\rangle  \langle s| A^*\}_{s \in \mathscr{S}}$ corresponds to the POVM on $H$ for embedded samples $\{A| s\rangle \}_{s \in \mathscr{S}}$. Given the density operator $\rho_{A}$, the probability density $p( A| s\rangle  |\rho_{A})$ can be derived using Born rule
\begin{equation}
\begin{aligned}
\label{eq:ker_3}
p \Big(A| s\rangle \big|\rho_{A} \Big) &= \mathrm{tr} \Big( A| s\rangle \langle s | A^* \, \rho_A \Big) \\
&= \frac{1}{\mathrm{tr}\left( A \, \rho \, A^* \right)} \, \int_{\mathscr{S}} \zeta(s') \, \big|\langle s| K |s'\rangle \big|^2  \, d\mu(s')
\end{aligned}
\end{equation}
This last equation is analogous to the usual kernel trick \citep[e.g.][]{epanechnikov1969non} and somehow relates to the kernel embedding of distribution \citep[e.g.][]{scholkopf2002learning}. The difference with the classic kernel embedding is the use of the square of the kernel $| \langle s'| K |s\rangle |$ rather than the the kernel itself. The reason lies in the nature of the wave functions which are probability amplitudes rather than probability densities. In equation \ref{eq:ker_3}, wave functions $\{|s\rangle\}_{s \in \mathscr{S}}$ associated to sample points are mapped on the image of $A$, whereas in kernel embedding Dirac measures $\{\delta_s\}_{s \in \mathscr{S}}$ associated to sample points are directly mapped on the image of $A$. As a result, the kernel is squared in equation \ref{eq:ker_3}. Embedding sample wave functions rather than densities guarantees that the $\rho_A$ can be interpreted as probability densities allowing to easily derive $p \Big(A| s\rangle \big|\rho_{A} \Big)$. While this difference may appear benign, it offers a natural framework which makes use of functional analysis while preserving the probabilistic nature of density operators.\newline
\\
\textbf{Homogeneous prior with non-noisy singleton.} \textit{It is straightforward to verify that for the density operator $\tilde{\rho}_{A} \propto A \, \tilde{\rho} \, A$ where $\tilde{\rho}$ is given by equation \ref{eq:res_2}, equation \ref{eq:ker_3} reduces to
\begin{equation}
\label{eq:ker_5}
\boxed{p\Big(A| s\rangle \big| \tilde{\rho}_{A}\Big) = \frac{1}{N \, \mathrm{tr}\left( A \, \rho \, A^* \right)} \, \sum_{i \in I}  \big| \langle S_i|K |s\rangle \big|^2}
\end{equation}
As with the usual kernel trick, it is not necessary to explicitly map samples on each of the wave functions in $\{|\psi_j \rangle \}_{j \in L}$ to compute $p( A| s\rangle  |\rho_{A})$. All which is needed is to be able to compute $| \langle S_i| K |s\rangle |^2$ for each singleton $\{S_i\}_{i \in I}$.}\newline
\\
Using equation \ref{eq:ker_3} it is now possible, and natural, to use Bayes' rule on the subset $\{\rho_A\}_{\rho \in \mathscr{E}}$ to find $P \left(\rho_A \big| \{A| S_i\rangle\}_{i \in I}  \right)$, i.e. the expression for Bayesian statistical learning of density operators $\rho_A$
\begin{equation}
\begin{aligned}
\label{eq:ker_4}
&P \left(\rho_A \big| \{A| S_i\rangle\}_{i \in I}  \right) \propto P \left( \rho_A  \right) \, \prod_{i \in I}  \int_{\mathscr{S}} P\left(S_i|s\right) \,  p \left(A| s\rangle \big| \rho_A \right) \, d\mu(s)
\end{aligned}
\end{equation}
where $P \left( \rho_A  \right)$ is the prior on the subspace of density operators $\rho_A$ and $P\left(S_i|s\right)$ describes the usual noise endowing each singleton as in equation \ref{eq:learning0}. Equation \ref{eq:ker_4} provides the general expression to learn density operators with embedded samples.\newline
\\
\textbf{Orthogonal projections.} \textit{We consider the special case where $A$ corresponds to an orthogonal projection operator $A$ such that,
\begin{equation}
\begin{aligned}
\label{eq:ker_5_1}
A = \sum_{j \in L} |\psi_j \rangle \langle \psi_j |
\end{aligned}
\end{equation}
where $L$ denotes a subset $L \in J$, then $A =A^*=A^* A$ and $\rho_A= A \, \rho_A \, A^*$. Using results from this section, it is possible to directly derive statistical learning equations of projection of density operators on subspaces spanned by $\{A \, \rho \, A^*\}_{\rho \in \mathscr{E}}$. Note that since coefficients $\{w(\psi_j,\psi_l)\}_{j,l \in L}$ are equal for both $\rho$ and $\rho_A$ then $\rho_A \to \rho$ as $L \to J$.}\newline
\\
Learning density operators using embedded samples allows to preserve the Bayesian formalism by embedding wave functions rather than probability densities. As with any kernel methods, using a subspace of all density operators reduces the complexity of the original Bayesian statistical learning problem at the expanse of discarding some information contained in the singletons. In order to illustrate these results, the next section details how embedding of density operators can be used to learn discrete orthogonal wavelet transform of density operators.

\subsection{Discrete orthogonal wavelet transform of density operators}

\begin{figure}
    \centering
    \begin{subfigure}[b]{0.47\textwidth}
        \includegraphics[width=\textwidth]{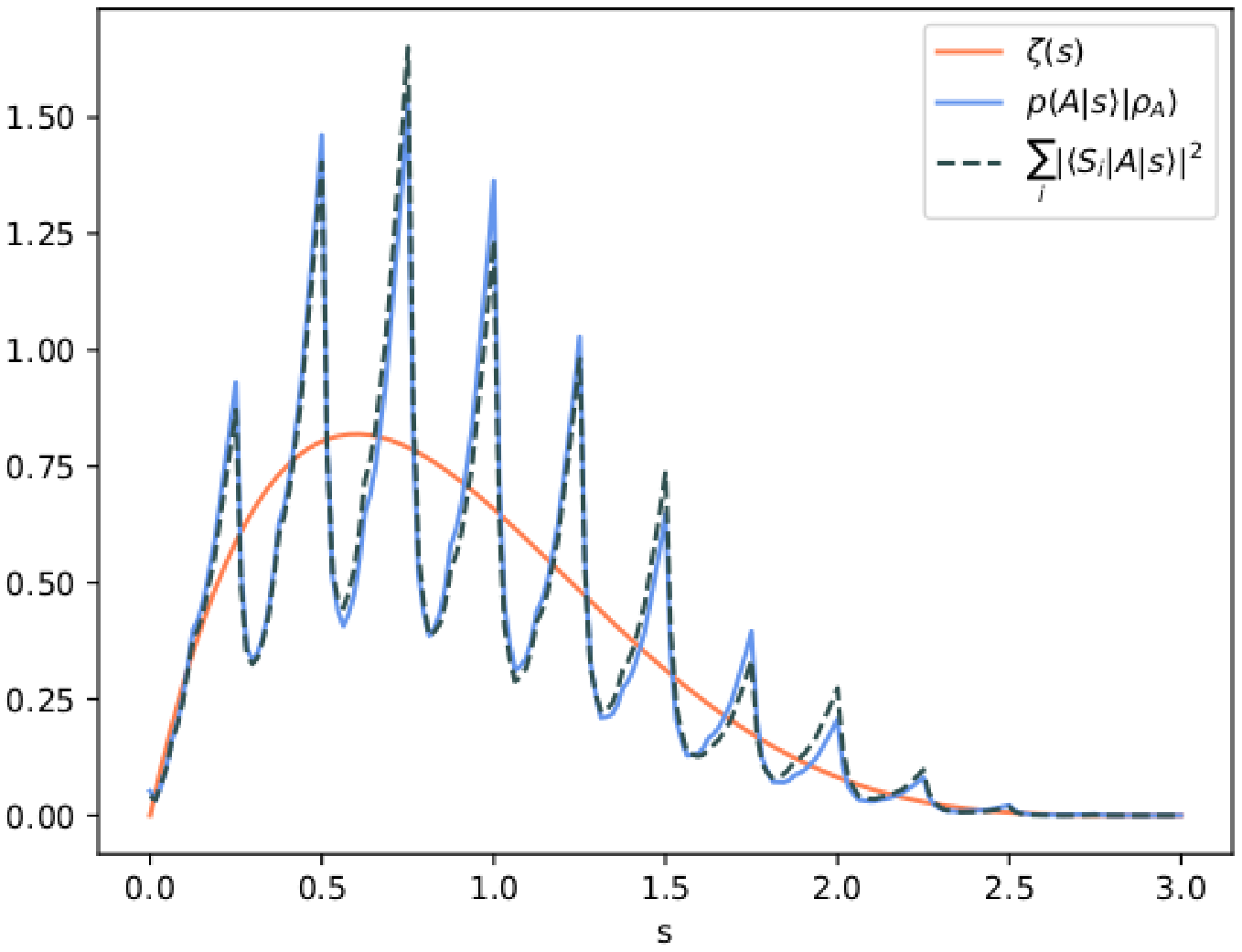}
        \caption{}
    \end{subfigure}
    ~ \quad 
    \begin{subfigure}[b]{0.47\textwidth}
        \includegraphics[width=\textwidth]{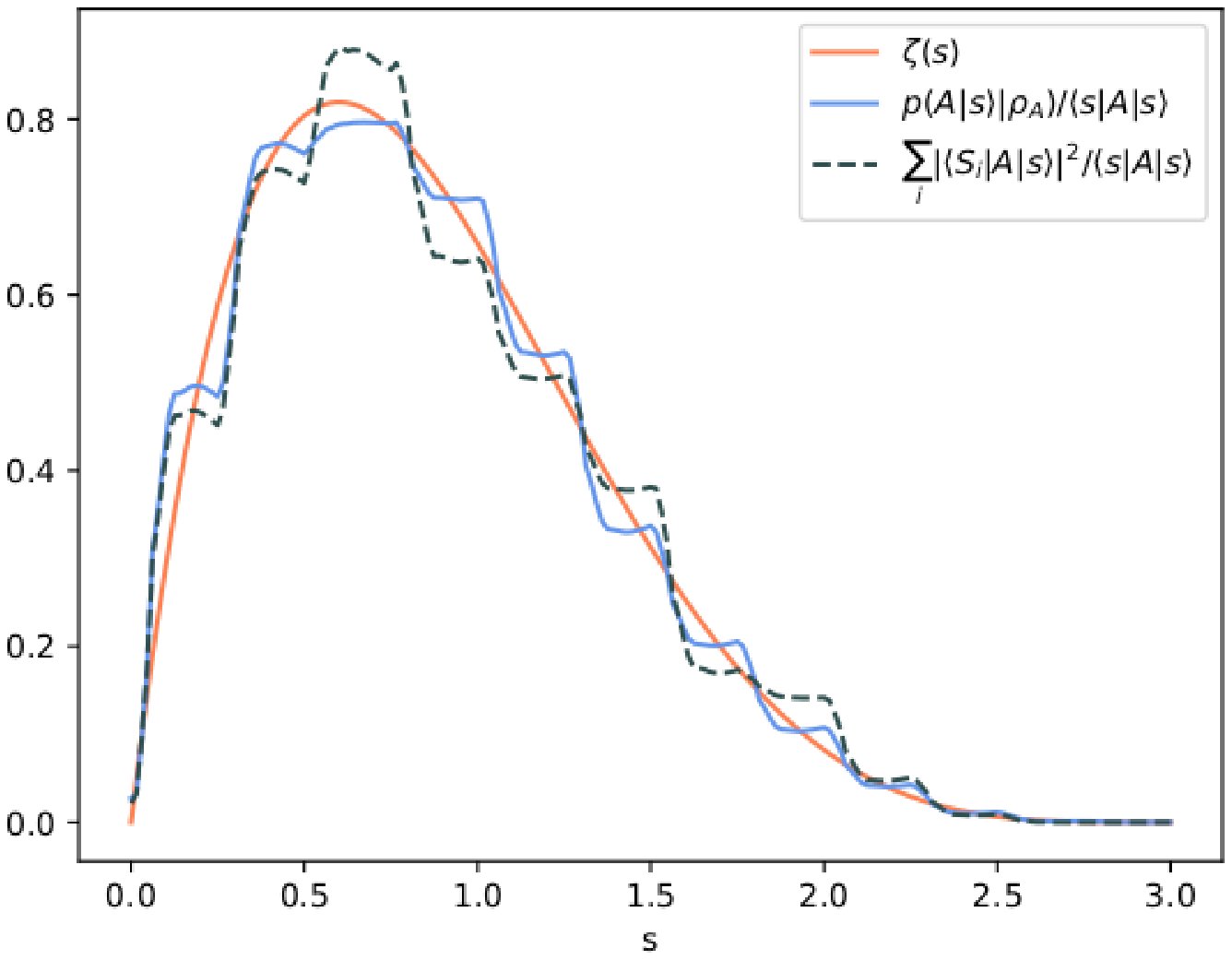}
        \caption{}
    \end{subfigure}
    \caption{(a) Target probability density $\zeta(s)$ with $p(A|s\rangle | \rho_A )$ and $p(A|s\rangle | \tilde{\rho}_A )$ where $p\left(A| s\rangle \big| \tilde{\rho}_{A}\right)$ is given by equation \ref{eq:ker_5}, (b) target probability density $\zeta(s)$ with $p(A|s\rangle | \rho_A )/\langle s| A|s\rangle$ and $p(A|s\rangle | \tilde{\rho}_A )/\langle s| A|s\rangle$.}\label{fig:ex2}
\end{figure}

Discrete orthogonal wavelet transform is a powerful tool which is extensively used in signal processing and machine learning \citep[e.g.][]{akansu2010emerging}. Its main advantage lies in the orthogonal property of the nested scale subspsaces which allows to sequentially improve the reconstruction of functions \citep[e.g.][]{Farge1992}. Consider the sample space $\mathscr{S}=\mathbb{R}$ with its associated Hilbert space $L^2(\mathbb{R},ds)$. For a given discrete orthogonal wavelet transform, the wavelet approximation at scale $2^{-n}$ of a function $f$ of $L^2(\mathbb{R},ds)$ is given by
\begin{equation}
\label{eq:w1}
\begin{aligned}
f_j(s) &= \sum_{k=-\infty}^{+\infty} \langle \phi_{nk} | f\rangle \, \phi_{nk}(s) 
\end{aligned}
\end{equation}
where $\phi_{nk}$ is the usual wavelet's father at scale $2^n$ with discrete translate $k$. By construction, $\{\phi_{nk}\}_{k}$ is a set of mutually orthonormal wave functions on $\mathscr{S}$. We introduce the following projection operator $A$
\begin{equation}
\label{eq:w2}
\begin{aligned}
A =  \sum_{k=-\infty}^{+\infty} |\phi_{nk} \rangle \langle \phi_{nk} |
\end{aligned}
\end{equation}
The projection operator $A$ is a linear map acting on $L^2(\mathbb{R},ds)$ which can be used to project density operators on the subspace induced by $\{\phi_{nk}\}_{k}$. 
\begin{figure}
    \centering
    \begin{subfigure}[b]{0.47\textwidth}
        \includegraphics[width=\textwidth]{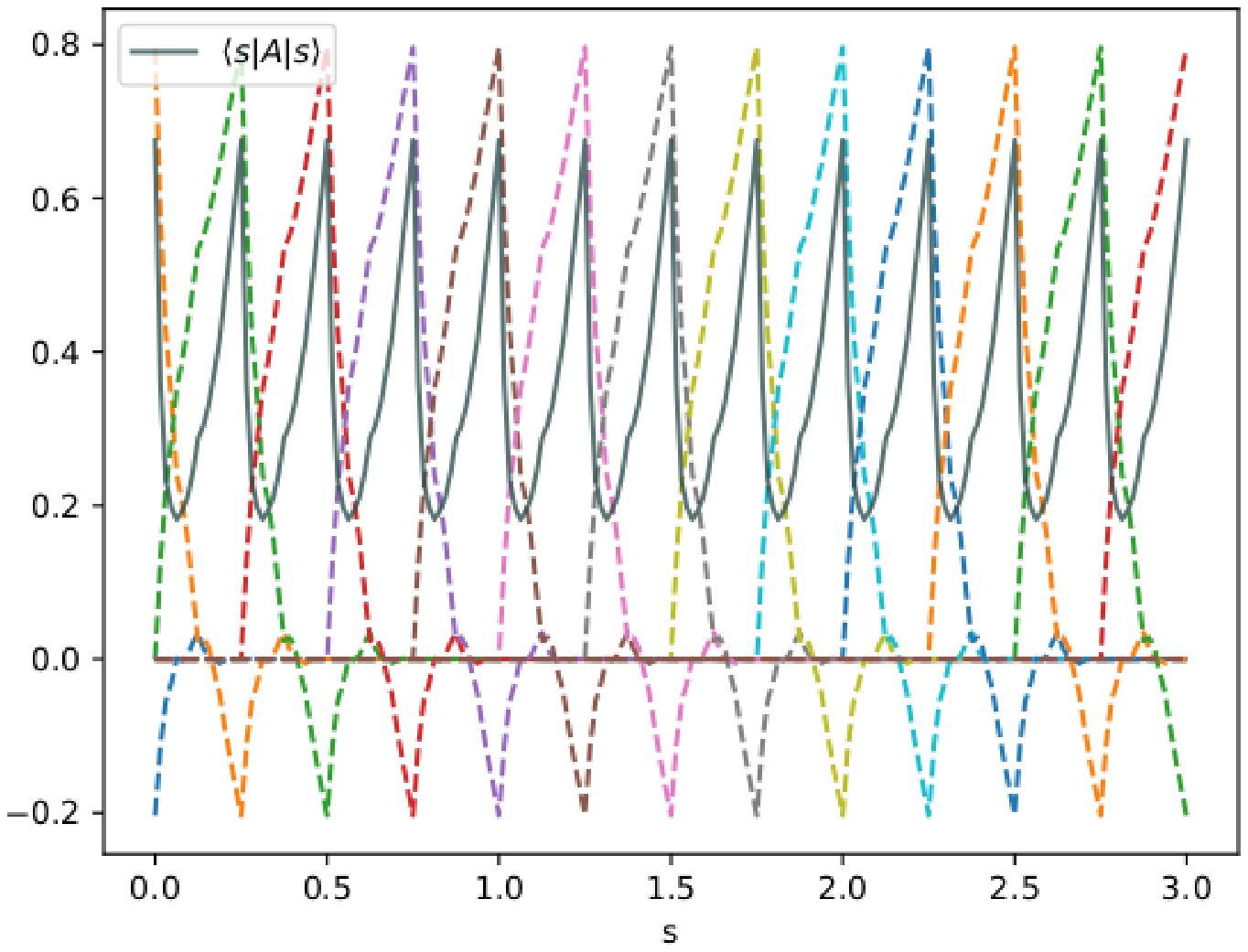}
        \caption{}
    \end{subfigure}
    ~ \quad 
    \begin{subfigure}[b]{0.47\textwidth}
        \includegraphics[width=\textwidth]{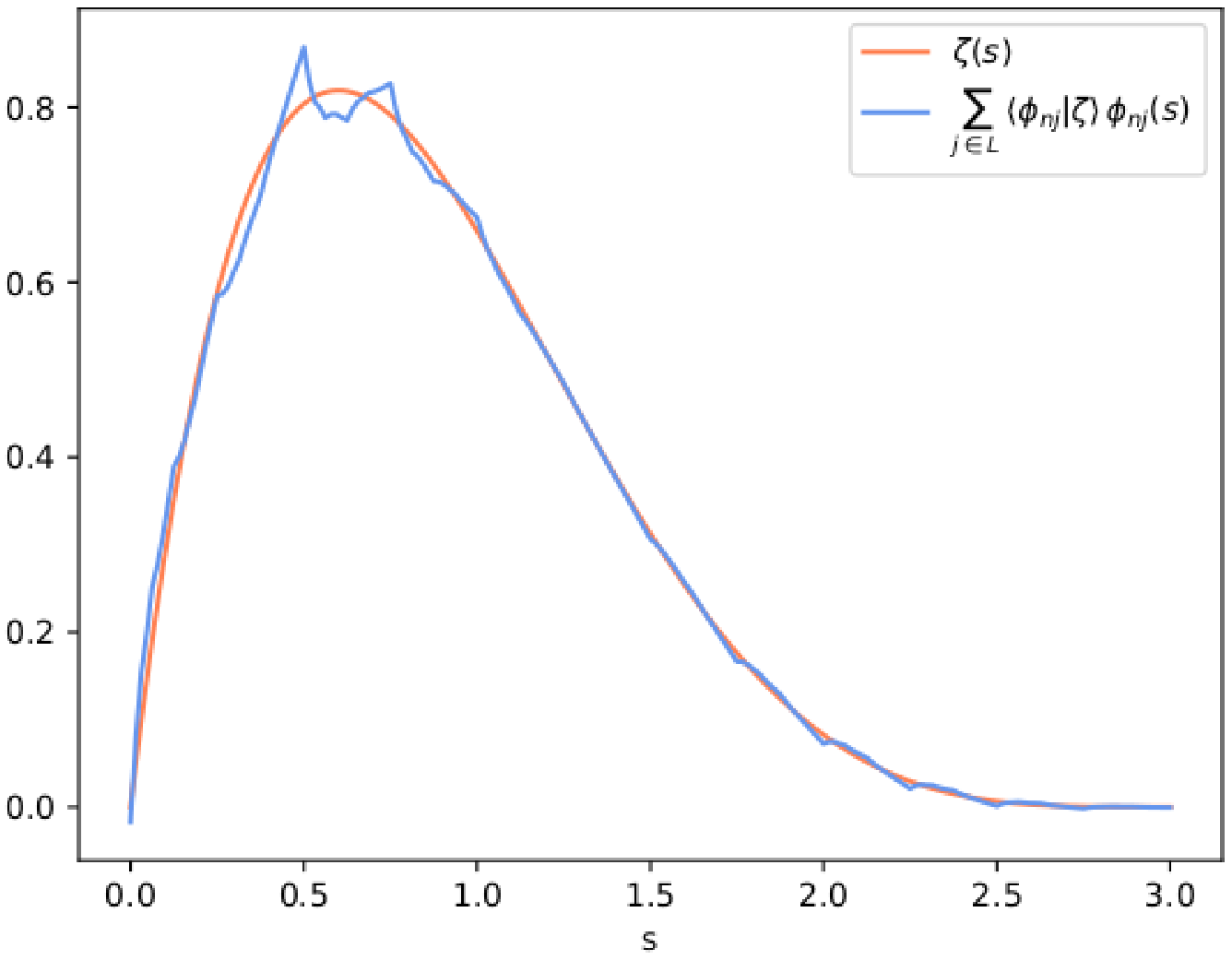}
        \caption{}
    \end{subfigure}
    \caption{(a) Normalized father wavelets and function $\langle s| A|s\rangle$, (b) probability density on $\mathscr{S}$ and its associated wavelet transform.}\label{fig:ex1}
\end{figure}
For a given sample $|s\rangle$ and a given density operator $\rho$, $A|s\rangle$ corresponds to the discrete orthogonal wavelet transform of $|s\rangle$ and $\rho_A$ corresponds to discrete orthogonal wavelet transform of $\rho$.\newline
\\
In order to illustrate the results, we choose the original sample space such that $\mathscr{S}=[0,3]$ and we use a beta probability density $\zeta(s)$. A set of 1D Daubechies tap 4 \citep[e.g.][]{daubechies1992ten} where the father's wavelet scale is chosen such that $n=2$ constitutes the orthogonal basis for the subspace of density operators. Associated father wavelets $\phi_{2k}$ are displayed in Figure \ref{fig:ex1}(a) along with $\mathrm{tr} \Big(A |s\rangle \langle s| A\Big)$. We voluntarily choose this family of non-symmetric wavelets to illustrate equation \ref{eq:ker_3} and the approximation $p(A|s\rangle | \rho_A )/ \langle s| A|s\rangle$ of the target probability density. Figure \ref{fig:ex1}(b) shows the target probability density $\zeta(s)$ along with its wavelet approximation at scale $2^{-n}$. A set of $N=300$ non-noisy independent singletons $\{S_i\}_{i \in I}$ was randomly generated from the target probability density $\zeta$. These non-noisy singletons were used to find the mode $\tilde{\rho}_{A}$ of the posterior probability distribution as given by equation \ref{eq:ker_5}. Probability densities $p(A|s\rangle | \rho_{A} )$ and $p(A|s\rangle | \tilde{\rho}_{A} )$ are shown in Figure \ref{fig:ex2}(a). It can be seen that $p(A|s\rangle | \rho_A )$ does not obviously relate to the wavelet approximation of the probability density. Indeed, when using density operators, singletons are treated as quantum particles described by their wave functions. Consequently, wave functions (rather than probability densities) are decomposed on the orthonormal basis and the kernel used corresponds to $\big|\langle s'| A |s\rangle \big|^2$ (rather than $\langle s'| A |s\rangle$). Using $\rho_A$ in place of $\rho$ corresponds to mapping all samples $\{|s\rangle\}_{s \in \mathscr{S}}$ to $\{A|s\rangle\}_{s \in \mathscr{S}}$. Note that according to equation \ref{eq:ker_2}, density operator $\rho_{A}$ can be understood as an ensemble of samples in state $\widehat{A| s\rangle}$ with probability $Z \left(s\right) \,\langle s |A| s\rangle$. In this last expression, $\langle s |A| s\rangle$ corresponds to the probability of observing samples in state $\widehat{A| s\rangle}$ assuming homogeneous probability density of generating samples in state $|s\rangle$. $\langle s |A| s\rangle$ varies greatly with $|s\rangle$ as can be seen in Figure \ref{fig:ex1}(a). We choose to display the normalized ratio of $p \left(A| s\rangle \big| \rho_A \right)$ over $\langle s |A| s\rangle$ in Figure \ref{fig:ex2}(b) to better highlight the link between $Z(s)$ and the density operator $\rho_{A}$.
This normalized ratio can be understood as the resulting approximation of $Z(s)$ induced by the mapping of samples $\{|s\rangle\}_{s \in \mathscr{S}}$ to $\{A|s\rangle\}_{s \in \mathscr{S}}$. While the projected density $\rho_A$ yields a proper probability density on the set of states $\{A|s\rangle\}_{s \in \mathscr{S}}$, the wavelet approximation of $\zeta(s)$ in Figure \ref{fig:ex1}(b) does not correspond to a probability density (it can take negative values). Decomposing wave functions rather than probability densities thus provides a natural framework which preserves the nature of density operators.

\section{Conclusion}

Using density operators representing ensembles of pure states of sample wave functions and wave functions in place of singletons provides a natural way to decompose samples on complete orthonormal basis of the sample space. Such decomposition allows to formulate the statistical Bayesian learning problem in different coordinate systems. Taking advantage of this representation, this study shows that it is possible to learn density operators on embedded sample spaces while preserving the Bayesian framework. It further proves that learning embedded density operators does not require to explicitly embed wave functions and can be performed efficiently using a kernel trick. One of the main advantage of decomposing wave functions rather than probability densities is the ability to preserve the nature of density operators throughout the process. The proposed approach thus differs from kernel embedding approaches which directly apply functional analysis tools to probability densities. This preliminary study remains largely incomplete and additional work is needed to figure out whether statistical Bayesian learning of density operators can result in novel efficient machine learning algorithms. In particular, applications to machine learning usual problems should be investigated and compared to state-of-the-art algorithms.

\bibliographystyle{apalike}
\bibliography{bib/references}

\begin{thebibliography}{}

\bibitem[Akansu et~al., 2010]{akansu2010emerging}
Akansu, A.~N., Serdijn, W.~A., and Selesnick, I.~W. (2010).
\newblock Emerging applications of wavelets: A review.
\newblock {\em Physical communication}, 3(1):1--18.

\bibitem[Barndorff-Nielsen et~al., 2003]{barndorff2003quantum}
Barndorff-Nielsen, O.~E., Gill, R.~D., and Jupp, P.~E. (2003).
\newblock On quantum statistical inference.
\newblock {\em Journal of the Royal Statistical Society: Series B (Statistical
  Methodology)}, 65(4):775--804.

\bibitem[Cabbolet, 2021]{cabbolet2021hyperreal}
Cabbolet, M.~J. (2021).
\newblock Hyperreal delta functions as a new general tool for modeling systems
  with infinitely high densities.
\newblock {\em Axioms}, 10(4):244.

\bibitem[Daubechies, 1992]{daubechies1992ten}
Daubechies, I. (1992).
\newblock {\em Ten lectures on wavelets}, volume~61.
\newblock Siam.

\bibitem[Epanechnikov, 1969]{epanechnikov1969non}
Epanechnikov, V.~A. (1969).
\newblock Non-parametric estimation of a multivariate probability density.
\newblock {\em Theory of Probability \& Its Applications}, 14(1):153--158.

\bibitem[{Farge}, 1992]{Farge1992}
{Farge}, M. (1992).
\newblock {Wavelet transforms and their applications to turbulence}.
\newblock {\em Annual Review of Fluid Mechanics}, 24:395--457.

\bibitem[Gonz{\'a}lez et~al., 2021]{gonzalez2021classification}
Gonz{\'a}lez, F.~A., Vargas-Calder{\'o}n, V., and Vinck-Posada, H. (2021).
\newblock Classification with quantum measurements.
\newblock {\em Journal of the Physical Society of Japan}, 90(4):044002.

\bibitem[Jarzyna and Ko{\l}ody{\'n}ski, 2020]{jarzyna2020geometric}
Jarzyna, M. and Ko{\l}ody{\'n}ski, J. (2020).
\newblock Geometric approach to quantum statistical inference.
\newblock {\em IEEE Journal on Selected Areas in Information Theory},
  1(2):367--386.

\bibitem[Li et~al., 2021]{li2021quantum}
Li, Q., Gkoumas, D., Lioma, C., and Melucci, M. (2021).
\newblock Quantum-inspired multimodal fusion for video sentiment analysis.
\newblock {\em Information Fusion}, 65:58--71.

\bibitem[Malley and Hornstein, 1993]{malley1993quantum}
Malley, J.~D. and Hornstein, J. (1993).
\newblock Quantum statistical inference.
\newblock {\em Statistical Science}, pages 433--457.

\bibitem[Paris, 2012]{paris2012modern}
Paris, M.~G. (2012).
\newblock The modern tools of quantum mechanics.
\newblock {\em The European Physical Journal Special Topics}, 203(1):61--86.

\bibitem[Phillips, 2013]{phillips2013introduction}
Phillips, A.~C. (2013).
\newblock {\em Introduction to quantum mechanics}.
\newblock John Wiley \& Sons.

\bibitem[Sch{\"o}lkopf et~al., 2002]{scholkopf2002learning}
Sch{\"o}lkopf, B., Smola, A.~J., Bach, F., et~al. (2002).
\newblock {\em Learning with kernels: support vector machines, regularization,
  optimization, and beyond}.
\newblock MIT press.

\bibitem[Sergioli et~al., 2018]{sergioli2018quantum}
Sergioli, G., Santucci, E., Didaci, L., Miszczak, J.~A., and Giuntini, R.
  (2018).
\newblock A quantum-inspired version of the nearest mean classifier.
\newblock {\em Soft Computing}, 22(3):691--705.

\bibitem[Tiwari and Melucci, 2019]{tiwari2019binary}
Tiwari, P. and Melucci, M. (2019).
\newblock Binary classifier inspired by quantum theory.
\newblock In {\em Proceedings of the AAAI conference on artificial
  intelligence}, volume~33, pages 10051--10052.

\bibitem[Wolf, 2006]{wolf2006learning}
Wolf, L. (2006).
\newblock Learning using the born rule.
\newblock Technical report.

\bibitem[Zhang et~al., 2018]{zhang2018quantum}
Zhang, Y., Song, D., Zhang, P., Wang, P., Li, J., Li, X., and Wang, B. (2018).
\newblock A quantum-inspired multimodal sentiment analysis framework.
\newblock {\em Theoretical Computer Science}, 752:21--40.

\end{thebibliography}

\end{document}